\documentclass[a4paper,12pt]{article}

\usepackage{amsmath,amsfonts}
\usepackage{mathrsfs}

\begin{document}

\def \Liminf{\mathop{\underline{\lim}}\limits}
\def \Limsup{\mathop{\overline{\lim}}\limits}
 \def\1{\mbox{1\hspace{-.25em}I}}
\def\Ex{{\bf E}}
\def\Pb{{\bf P}}
\let\bb\mathbb 
\def\UU{{\bb U}}
\def\KK{{\bb K}}
\def\BB{{\bb B}}
\def\CC{{\bb C}}
\def\NN{{\bb N}}
\def\AA{{\bb A}}
\def\DD{{\bb D}}
\def\RR{{\cal R}}
\def\sgn{{\rm sgn}}
\let\scr\mathscr
\newtheorem{theorem}{Theorem}
\newtheorem{lemma}{Lemma}
\newtheorem{proposition}{Proposition}

\title{{ On frequency estimation of periodic  ergodic diffusion process}}
\author{R. \textsc{H\"opfner} \\
{\small Johannes Gutenberg Universit\"at Mainz}\\[8pt]
Yu. A. \textsc{Kutoyants}\\
{\small  Universit\'e du Maine, Le Mans}}
\date{ }
\maketitle

\begin{abstract}
We consider the problem of frequency estimation by observations of the
periodic diffusion process possesing ergodic properties in two different
situations. The first one corresponds to continuously diffeentiable with
respect to parameter  trend
coefficient and the second - to discontinuous trend coefficient. It is shown
that in the first
case the maximum likelihood and bayesian estimators are asymptotically normal
with rate $T^{3/2}$ and in the second case these estimators have different
limit distributions with the rate $T^2$. 
\end{abstract}



\section{Introduction}

Let us consider the model ``signal in noise'' of the following type
\begin{equation}
\label{1}
x\left(t\right)=S\left(\vartheta ,t\right)+n\left(t\right),\quad 0\leq t\leq T,
\end{equation}
where $S\left(\vartheta ,\cdot \right)$ is the signal transmitting the
``information'' $\vartheta $ and observed in the presence of additive
``noise'' $n\left(\cdot \right)$. This is  a typical model for the theory of
telecommunications.  There is a  alarge diversity of statistical
problems related to  this model. One way is to study the different noises
(white Gaussian, ``colored'' Gaussian, stationary {\it etc.}) and another way
is to study the different types of ``modulations'', i.e.; to choose  the signals
$S\left(\cdot \right)$, like {\it amplitude modulation} $S\left(\vartheta,t\right)=\vartheta 
h\left(t\right)$, {\it phase modulation} $S\left(\vartheta,t
\right)=h\left(t-\vartheta \right)$ or {\it frequency modulation} $S\left(\vartheta,t
\right)=h\left(\vartheta t \right)$. Here the function $h\left(\cdot \right)$
is usualy supposed to be periodic. The most developed is the theory of
estimation of periodic signals observed in white Gaussian noise.

The problem of period (or frequency) estimation has particularities which put
it in some sense out of traditional ($\sqrt{n}$) statistical framework. Let us
remind some known properties of the maximum likelihood estimator (MLE) of the
frequency $\vartheta \in \left(\alpha ,\beta \right),0<\alpha <\beta <\infty $
obtained by Ibragimov and Khasminskii \cite{IH81} for the model {\it signal in
  white Gaussian noise} (SWN) and some related problems for inhomogeneous Poisson
processes.

Suppose that the observed process is 
$$
 {\rm d}X_t=A\sin\left(\vartheta t\right){\rm d}t+\sigma {\rm d}W_t,\quad
X_0=0,\quad 0\leq t\leq T
$$ 
(SWN) and we have to estimate the frequency $\vartheta $ by the observations
$X^T=\left(X_t,0\leq t\leq T\right)$. We are interested by the properties of
estimators in asymptotics of ``large samples'': $T\rightarrow \infty $. 
 The Fisher information is
$$
{\rm I}_T\left(\vartheta \right)=\frac{A^2}{\sigma
  ^2}\int_{0}^{T}t^2\cos^2\left(\vartheta t\right)\,{\rm d}t=\frac{A^2T^3}{3\sigma
  ^2} \,\left(1+o\left(1\right)\right)
$$
and the MLE $\hat\vartheta _T$ is asymptotically normal with the rate $T^{3/2}$, i.e.;
\begin{equation}
\label{2}
 T^{3/2}\left(\hat \vartheta _T-\vartheta \right)\Rightarrow {\cal
  N}\left(0,\frac{3\sigma ^2}{A^2}\right),\qquad  \Ex_\vartheta \left(\hat 
\vartheta _T- \vartheta\right)^2=\frac{3\sigma
  ^2}{A^2T^3}\,\left(1+o\left(1\right)\right).  
\end{equation}
Note that if $\beta =\infty $, then the (uniformly in $\vartheta $) consistent
estimation is impossible.  Even if we allow $\beta _T\rightarrow \infty $,
then for $\beta _T<\exp\left\{\left(\frac{A^2}{4\sigma ^2}-\varepsilon
\right)T\right\}$ (any $\varepsilon >0$) the MLE is consistent and for $\beta
_T>\exp\left\{\left(\frac{A^2}{4\sigma ^2}+\varepsilon \right)T\right\}$ the
uniformly consistent estimation of $\vartheta $ is impossible (see \cite{IH81},
Section 7.1 for exact statemts and proofs).

If we consider the problem of parameter  estimation by
observations 
$$
{\rm d}X_t=S\left( t-\vartheta\right){\rm d}t+\sigma {\rm d}W_t,\quad
X_0=0,\quad 0\leq t\leq T,
$$ 
where $S\left(\cdot \right)$ is periodic function of period 1 having a discontinuity at
points $\tau_* +k, k=0,1,2,\ldots $, then the rate of convergence of the MLE
$\hat \vartheta _T$ is different. Let us denote $S\left(\tau_* -\right)$ and
$S\left(\tau_* +\right)$ the left and right limits, $S\left(\tau_*
+\right)-S\left(\tau_* -\right)=r\not=0$. Then 
\begin{equation}
\label{3}
T\left(\hat \vartheta _T-\vartheta \right)\Longrightarrow\eta ,\qquad \Ex_\vartheta
\left(\hat\vartheta _T-\vartheta \right)^2= \frac{26\,\sigma
  ^2}{r^2\,T^2} \,\left(1+o\left(1\right)\right), 
\end{equation}
 where $\eta $ is a random variable (see \cite{IH81}, Section 7.2).

  The
similar problems of parameter estimation were considered in \cite{Kut98} for
the model of periodic Poisson process. It was supposed that the observed
inhomogeneous Poisson process $X^T=\left(X_t,0\leq t\leq T\right)$ has
intensity function
$$
S\left(\vartheta ,t\right)=S\left(\vartheta t\right)
$$
where $S\left(t\right)$ is $\tau $-periodic smooth function.  It was shown
that the MLE $\hat \vartheta _T$ is asymptotically normal with 
the  rate $T^{3/2}$  :
\begin{equation}
\label{4}
 T^{3/2}\left(\hat \vartheta _T-\vartheta \right)\Longrightarrow {\cal
  N}\left(0,a^2\right),\qquad \Ex_\vartheta \left(\hat \vartheta
   _T-\vartheta\right)^2=\frac{a^2}{T^3} \,\left(1+o\left(1\right)\right),
\end{equation}
see \cite{Kut98}, Section 2.3 for details. 

In the case of discontinuous periodic intensity $S\left(t-\vartheta \right)$
(shift parameter) the rate is (like \eqref{3}) $T$, i.e.;
\begin{equation}
\label{5}
T\left(\hat \vartheta _T-\vartheta \right)\Longrightarrow \xi ,\qquad \Ex_\vartheta
\left(\hat\vartheta _T-\vartheta
\right)^2=\frac{c^2}{T^2}\left(1+o\left(1\right)\right).
\end{equation}
See \cite{Kut98}, Section 5.1 for details. Moreover, the problem of
frequency estimation of periodic discontinuous intensity function
$S\left(\vartheta t\right)$ was also considered and it was shown that the rate of
convergence of the MLE is $T^2$, i.e.; we have the limits

\begin{equation}
\label{6}
T^2\left(\hat \vartheta _T-\vartheta \right)\Longrightarrow \zeta ,\qquad \Ex_\vartheta
\left(\hat\vartheta _T-\vartheta
\right)^2=\frac{b^2}{T^4}\,\left(1+o\left(1\right)\right) .
\end{equation}

In the present work we consider the problem of frequency estimation in the
case of periodic discontinuous trend coefficient of ergodic diffusion
process. This work is a continuation of our study of parameter estimation problems for
periodic diffusion processes  started in \cite{HK10a}-\cite{HK10c}.
In all these works we suppose that the observed diffusion process is given by the equation 
\begin{equation}
\label{7}
{\rm d}X _t=[S\left(\vartheta ,{t}\right)+b\left(X
  _t\right)]\,
{\rm d}t+\sigma \left(X _t\right) \;{\rm d}W_t,\quad  X_0,\quad 0\leq t\leq T, 
\end{equation}
 where the function $S\left(\vartheta ,t\right),t\geq 0$ is (known) periodic
 of period $\tau $, i.e., $S\left(\vartheta ,t+k\tau \right)=S\left(\vartheta
 ,t\right) $ the functions $b\left(\cdot \right)$ and $\sigma \left(\cdot
 \right)$ are known and smooth.

This problem of parameter estimation can be considered as particular
 case of \eqref{1} with ``diffusion noise''
 $n\left(t\right)=b\left(X_t\right)+\sigma \left(X_t\right)\dot
 W_t$. Therefore once more we have a problem of the theory of
 telecommunication (transmission of signals) but this model of observations
 is as well interesting in some biological experiments related to membrane
 potential data sets (see H\"opfner \cite{H}).  We suppose that the diffusion
 process has ergodic properties and describe the asymptotics of the
 MLE and BE in regular and singular (discontinuous)  situations. The existance
 of periodic solution for Markov processes with periodic coefficients were
 studied by Khasminskii \cite{H} and the ergodic properties (law of large
 numbers, periodic invariant density...) used in the present
 work were  obtained by H\"opfner and L\"ocherbach \cite{HL10}. 

We have to note that if the diffusion coefficient is a deterministic function,
say, $\sigma \left(x\right)\equiv \sigma >0$, then the simple transformation
(we suppose always that $b\left(x\right)$ is known)
$$
Y_t=X_0-\int_{0}^{t}b\left(X_s\right)\,{\rm d}s
$$
reduces the equation \eqref{7} to the well known signal in WGN model
$$
{\rm d}Y _t=S\left(\vartheta ,{t}\right)\,
{\rm d}t+\sigma  \;{\rm d}W_t,\quad  Y_0=0,\quad 0\leq t\leq T, 
$$
and for this model all mentioned above problems are already well studied.

Let us denote by $\left\{{\bf P}_\vartheta ^{\left(T\right)},\vartheta \in \Theta \right\}$ the
family of measures induced by the solutions of \eqref{7} in the measurable space $\left({\cal
  C}\left[0,T\right],{\scr B}\left[0,T\right]\right)$ and put
$$
L\left(\vartheta ,X^T\right)=\frac{{\rm d}{\bf P}_\vartheta
  ^{\left(T\right)}}{{\rm d}{\bf P} ^{\left(T\right)}}\left(X^T\right) ,
$$
where ${\bf P} ^{\left(T\right)} $ is the measure corresponding to the process
\eqref{7} with $S\left(\vartheta ,t\right)\equiv 0$. 
The likelihood ratio is
\begin{equation}
\label{8}
L\left(\vartheta ,X^T\right)=\exp\left\{\int_{0}^{T}\frac{S\left(\vartheta
  ,t\right)}{\sigma \left(X_t\right)^2} {\rm d}X_t-\int_{0}^{T}\frac{S\left(\vartheta
  ,t\right)^2+2S\left(\vartheta, t\right)b\left(X_t\right)}{2\sigma
  \left(X_t\right)^2} {\rm d}t\right\} 
\end{equation}
and the estimators are defined by the usual formulas: for the MLE
$\hat\vartheta _T$ we have
\begin{equation}
\label{MLE}
L\left(\hat\vartheta _T,X^T\right)=\sup_{\theta \in \Theta }L\left(\theta,X^T\right),
\end{equation}
and for Bayesian estimator $\tilde\vartheta _T$ we write
\begin{equation}
\label{BE}
\tilde\vartheta _T=\frac{\int_{\Theta }^{}\theta p\left(\theta
  \right)L\left(\theta ,X^T\right){\rm d}\theta }{\int_{\Theta }^{}
  p\left(\theta \right)L\left(\theta ,X^T\right){\rm d}\theta} ,
\end{equation}
i.e., we suppose that the loss function is quadratic and the density a pripory
$p\left(\cdot \right)$ is given. We study the asymptotic properties of these
estimators with the help of the methode by Ibragimov and Khasminskii
\cite{IH81} which consists in the establishing some properties of the
normalized likelihod ratio process
$$
Z_T\left(u\right)=\frac{L\left(\vartheta +\varphi
  _Tu,X^T\right)}{L\left(\vartheta ,X^T\right)},\qquad u\in
\UU_T=\left[\frac{\alpha -\vartheta }{\varphi _T}, \frac{\beta  -\vartheta }{\varphi _T}\right]
$$
where $\vartheta $ is the true value and the choice of the normalizing
function $\varphi _T\rightarrow 0$ depends on the ``smoothness'' of the
problem. 

In the first work \cite{HK10a} we supposed that the function $S\left(\vartheta
,t\right)=S\left(t-\vartheta \right)$ and $S\left(t\right),t\geq 0$ is
periodic function having a jump $r=S\left(\tau_*
+\right)-S\left(\tau_*-\right)\neq 0$ at the point $\tau _*\in \left(0,\tau
\right)$. It is shown that the choice of the function $\varphi _T=T^{-1} $
provides the limit
\begin{equation}
\label{9}
Z_T\left(u\right)\Longrightarrow Z\left(u\right)=\exp\left\{\gamma
W\left(u\right)-\frac{\left|u\right|}{2}\gamma ^2\right\} ,\qquad u\in \RR
\end{equation}
with some constant $\gamma $. Here $W\left(u\right)$ is double sided Wiener
process.
 The estimators have the following properties. Let us put
\begin{equation}
\label{10}
Z\left(\hat u\right)=\sup_u Z\left(u\right),\qquad \tilde
u=\frac{\int_{-\infty }^{\infty }u\,Z\left(u\right)\,{\rm d}u }
{\int_{-\infty }^{\infty }Z\left(u\right)\,{\rm d}u},
\end{equation}
then we can write
\begin{equation}
\label{11}
T\left(\hat\vartheta _T-\vartheta \right)\Longrightarrow \hat u,\qquad
T\left(\tilde\vartheta _T-\vartheta \right)\Longrightarrow \tilde u, 
\end{equation}
and convergence of all polynomial moments of these estimators (like \eqref{3})
take place (see \cite{HK10a}).

In the second work \cite{HK10b} we considered the usual (regular) estimation problem
with smooth periodic function $S\left(\vartheta ,t\right)$ such that its derivative is
periodic too. It is shown that with classical normalization $\varphi _T=T^{-1/2}$ the
corresponding family of  measures is LAN :
\begin{equation}
\label{12}
Z_T\left(u\right)\Longrightarrow Z\left(u\right)=\exp\left\{u\Delta
-\frac{u^2}{2}{\rm I}\left(\vartheta \right)\right\},
\end{equation}
where ${\rm I}\left(\vartheta \right) $ is the Fisher information (on one
period) and $\Delta \sim {\cal N}\left(0,{\rm I}\left(\vartheta
\right)\right)$. For the estimators we obtain, as usual, asymptotic normality
\begin{equation}
\label{13}
\sqrt{T}\left(\hat\vartheta _T-\vartheta \right)\Longrightarrow \frac{\Delta
}{{\rm I}\left(\vartheta \right)},\qquad \sqrt{T}\left(\tilde\vartheta
_T-\vartheta \right)\Longrightarrow \frac{\Delta }{{\rm I}\left(\vartheta
  \right)}
\end{equation}
and convergence of all  polynomial moments
$$
\Ex_\vartheta \left(\hat\vartheta _T-\vartheta \right)^2=\frac{1}{T{\rm I}\left(\vartheta
  \right) }\left(1+o\left(1\right)\right).
$$
The last work is devoted to the study of the local structure of the family of
measures corresponding to the observations \eqref{7} with
$S\left(\vartheta,{t}\right)=S\left(\vartheta{t}\right) $ 
where $S\left(t\right)$ is a periodic function. We
describe the asympptotic behavior of the normalized likelihood ratio
$Z_T\left(u\right)$ in two situations: when the periodic function
$S\left(t\right)$ is smooth and when it is discontinuous. It is shown that in
the first case the normalizing function $\varphi _T=T^{-3/2} $ and the limit
is like \eqref{12} and in the second case $\varphi _T=T^{-2} $  and the limit
is like \eqref{9}. In the present work we descride the properties of the MLE
and Bayesian estimators of the frequency  in the same two situations. We show
that in the smooth case the MLE and BE are asymptotically normal similar to
\eqref{4} and in the discontinuous case  we have convergence like \eqref{6}
but with the different limit distributions.

\section{Main results}

The observed diffusion process$X^T=\left(X_t,0\leq t\leq T\right)$ satisfies
the stochastic differential equation 
\begin{equation}
\label{14}
{\rm d}X _t=[S\left(\vartheta{t}\right)+b\left(X
  _t\right)]\,
{\rm d}t+\sigma \left(X _t\right) \;{\rm d}W_t,\quad  X_0,\quad 0\leq t\leq T, 
\end{equation}
and we study the properties of estimators of the  parameter $\vartheta $ in two
situations : when the function $S\left(t\right)$ is smooth (regular estimation
problem) and when the function $S\left(t\right)$ has a discontinuity (singular
estimation problem). 

In both cases we suppose that the following conditions are  fulfilled.
\bigskip

{\bf A1.} {\it The function  $S\left(t\right)$ is  bounded and  periodic  of period $\tau
=1$. The functions $b\left(\cdot \right),\sigma \left(\cdot \right)\in {\cal
  C}_b^3$, i.e.; have three continuous bounded derivatives. 
The parameter $\vartheta \in \left(\alpha ,\beta \right)=\Theta $}.

Let us introduce the constants $m, M$ by the relation 
$$
m\leq S\left(t\right)\leq M,\quad \quad t\in \left[0,1
  \right],
$$
and  put $S_{-}=\min \left(m,0\right)$ and $S_+=\max\left(M,0\right)$.

 {\bf A2.} {\it  There exist constants $A>0$  and $\varepsilon >0$
such that for $\left|x\right|>A$ we have }
$$
2x\1_{\left\{x<-A\right\}}S_{-} +2x\1_{\left\{x>A\right\}}S_{+}
+2xb\left(x\right)+\sigma \left(x\right)^2<-\varepsilon  .
$$

Note that by condition {\bf A1} these
functions satisfy  the global Lipshitz condition 
$$
\left|b\left(x \right)-b\left(y\right)\right|+\left| \sigma
\left(x\right)-\sigma \left(y \right)\right| \leq L\;\left|x-y\right|
$$
and  the linear growth  condition  
$$ \left|b\left(x\right)\right|+\left|\sigma \left(x\right)\right|\leq
L_1\left(1+\left|x\right|\right).
$$
Therefore the equation \eqref{14} has a unique strong solution \cite{LS}.

 {\bf A3.} {\it The  diffusion coefficient  is a bounded function separated from
zero : there exist two constants $k,K$ such that }
\begin{equation}
\label{15}
 0<\kappa \leq \sigma \left(x\right)^2\leq K  
\end{equation}

Under  conditions {\bf A1},{\bf A2},{\bf A3}  the diffusion process
has ergodic properties (H\"opfner 
and L\"ocherbach, \cite{HL10} ), i.e., there exists an invariant (periodic) density
function $f_\vartheta \left(t,x\right)$ such that for any absolutely
integrable $\tau =1/\vartheta $-periodic in time  function $h\left(\vartheta ,t,x\right)$ we have (with probability 1) the
following limits
\begin{align}
\label{8a}
&\frac{1}{T}\int_{0}^{T}h\left(\vartheta ,t,X _t\right)\,{\rm d}t\longrightarrow\frac{1}{\tau }
\int_{-\infty }^{\infty }\int_{0}^{\tau }h\left(\vartheta ,t,x\right)
\;f_\vartheta \left(t,x\right){\rm d}t\;{\rm d}x, \\
&
\frac{1}{n}\sum_{k=1}^{n}h\left(\vartheta ,t_*+k\tau ,X_{t_*+k\tau }\right)\longrightarrow
\int_{-\infty }^{\infty }h\left(\vartheta ,t_*,x\right)f_\vartheta \left(t_*,x\right){\rm d}x.
\label{8b}
\end{align}

To prove the asymptotic efficiency we need the following uniform law of large
numbers.

{\bf A4.} {\it The convergence in \eqref{8a}, \eqref{8} is uniform on compacts $\KK\in
\Theta $.} 
\bigskip

A sufficient for {\bf A4} condition is given in the Section 3.
Below the condition {\bf A}= ({\bf A1},{\bf A2},{\bf A3},{\bf A4}).

\subsection{Smooth trend}

We consider the problem of frequency $\vartheta $ estimation by observations $X^T$
of the periodic diffusion process \eqref{14}.

{\bf B.} {\it The periodic function $S\left(t \right),t\geq 0$ is  nonconstant
  and  continuously
  differentiable.}

 The role of Fisher information in our problem plays the quantity 
$$
{\rm I}\left(\vartheta\right)=\frac{1}{3\tau }\int_{0}^{\tau }\dot S\left(\vartheta t\right)^2  \int_{-\infty }^{\infty
} \frac{f_\vartheta\left({t},x\right) }{\sigma
  \left(x\right)^2}\;{\rm d}x \;{\rm d}t,
$$
where dot means derivation : $\dot S\left(t\right)={\rm d}S\left(t\right)/{\rm
  d}t$. 
 
Introduce the lower bound on the meansquare risk of all estimators
$\bar\vartheta _T$:
\begin{equation}
\label{bound1}
\Liminf_{\delta \rightarrow 0}\Liminf_{T \rightarrow \infty
}\sup_{\left|\vartheta -\vartheta _0\right|\leq \delta }T^3 \Ex_\vartheta
\left|\bar\vartheta _T-\vartheta \right|^2\geq {\rm
  I}\left(\vartheta_0\right)^{-1}.
\end{equation}
This is a version of the  well-known Hajek-Le Cam lower bound (see for example
\cite{IH81}). We call an
estimator $\vartheta _T^\star$  asymptotically efficient if for all $\vartheta
_0\in \Theta $ we have the equality
\begin{equation}
\label{asef}
\lim_{\delta \rightarrow 0}\lim_{T \rightarrow \infty
}\sup_{\left|\vartheta -\vartheta _0\right|\leq \delta }T^3 \Ex_\vartheta
\left|\vartheta _T^\star-\vartheta \right|^2= {\rm
  I}\left(\vartheta_0\right)^{-1}.
\end{equation}

\begin{theorem}
\label{T1} 
 Let the conditions {\bf A} and  {\bf B} be fulfilled. Then the  MLE $\hat \vartheta _T$ and
BE $\tilde\vartheta _T$ have the following properties uniformly on compacts
$\KK\subset\Theta $.
\begin{itemize}
\item These estimators are   consistent:   for any $\delta >0$
$$
\sup_{\vartheta \in \KK}\Pb_\vartheta \left\{\left|\hat\vartheta _T-\vartheta
\right|>\delta \right\}\rightarrow 0, \qquad \sup_{\vartheta \in
  \KK}\Pb_\vartheta \left\{\left|\tilde\vartheta _T-\vartheta 
\right|>\delta \right\}\rightarrow 0. 
$$
\item These estimators are asymptotically normal
\begin{equation*}
\label{17}
T^{3/2}\left(\hat\vartheta _T-\vartheta \right)\Longrightarrow \zeta ,\quad
T^{3/2}\left(\tilde\vartheta _T-\vartheta \right)\Longrightarrow \zeta ,\qquad
\zeta \sim {\cal 
  N}\left(0,{\rm I}\left(\vartheta \right)^{-1}\right).
\end{equation*} 
\item We have  the convergence  of  moments : for any $p>0$ 
$$
\lim_{T\rightarrow \infty }T^{\frac{3p}{2}}\;\Ex_\vartheta \left|\hat\vartheta
_T-\vartheta \right|^p=\Ex_\vartheta \left|\zeta \right|^p,\quad  \lim_{T\rightarrow
  \infty }T^{\frac{3p}{2}}\;\Ex_\vartheta \left|\hat\vartheta 
_T-\vartheta \right|^p=\Ex_\vartheta \left|\zeta \right|^p,
$$

\item  The both
estimators are asymptotically efficient in the sense \eqref{asef} .
\end{itemize}
\end{theorem}

{\bf Proof.} Let us introduce the normalized likelihood ratio
\begin{align*}
Z_T\left(u\right)=\frac{L\left(\vartheta
  +T^{-3/2}u,X^T\right)}{L\left(\vartheta ,X^T\right)},\qquad u\in
\UU_T=\left(T^{3/2}\left(\alpha -\vartheta 
\right),T^{3/2}\left(\beta  -\vartheta \right) \right).  
\end{align*} 
According to \eqref{8} it has the form (below $\vartheta _u=\vartheta +T^{-3/2}u$)
\begin{align*}
\ln Z_T\left(u\right)=\int_{0}^{T}\frac{S\left(\vartheta
  _ut\right)-S\left(\vartheta t\right)}{\sigma \left(X_t\right)}{\rm
  d}W_t-\frac{1}{2}\int_{0}^{T}\left(\frac{S\left(\vartheta
  _ut\right)-S\left(\vartheta t\right)}{\sigma \left(X_t\right)}\right)^2{\rm
  d}t.
\end{align*}
This proces can be written as 
\begin{align*}
\ln Z_T\left(u\right)=\frac{u}{T^{3/2}}\int_{0}^{T}\frac{t\,\dot S\left(\vartheta
  t\right)}{\sigma \left(X_t\right)}\;{\rm
  d}W_t-\frac{u^2}{2T^{3}}\int_{0}^{T}\left(\frac{t\,\dot S\left(\vartheta
  t\right)}{\sigma \left(X_t\right)}\right)^2{\rm d}t +r_T.
\end{align*}
It was shown (see \cite{HK10c}) that 
$$
\frac{1}{T^{3}}\int_{0}^{T}\left(\frac{t\,\dot S\left(\vartheta
  t\right)}{\sigma \left(X_t\right)}\right)^2{\rm d}t\longrightarrow {\rm
  I}\left(\vartheta \right) ,\qquad r_T\rightarrow 0
$$
and
$$
\Delta _T\left(\vartheta \right)=\frac{1}{T^{3/2}}\int_{0}^{T}\frac{t\,\dot S\left(\vartheta
  t\right)}{\sigma \left(X_t\right)}\;{\rm   d}W_t\Longrightarrow \Delta\left(\vartheta \right) \sim {\cal
  N}\left(0,{\rm I}\left(\vartheta \right)\right).
$$

 Moreover, as we suppose uniform law of large numbers \eqref{8}, this
convergence is uniform on compacts $\KK\subset \Theta $.  Therefore, if we
introduce the random process (see \eqref{12})
$$
Z\left(u\right)=\exp\left\{u\,\Delta \left(\vartheta
\right)-\frac{u^2}{2}{\rm I}\left(\vartheta \right)\right\},\qquad u\in \RR,
$$
then the following result take place.
\begin{lemma}
\label{L1} The finite dimensional distributions of the random process
$Z_T\left(\cdot \right)$ converge to the finite dimensional distributions of the
process $Z\left(\cdot \right)$ uniformly in $\vartheta \in \KK$.
\end{lemma}

For the proof see \cite{HK10c}, Theorem 1.1. Just note that in \cite{HK10c} we
do not supposed the uniformity of this convergence and at present it follows
from the uniform law of large numbers.

\begin{lemma}
\label{L2} For any $R>0$ the following inequaulity holds
\begin{equation}
\label{18}
\sup_{\vartheta \in \KK} \sup_{\left|u_1\right|+\left|u_2\right|\leq R}   \left|u_2-u_1\right|^{-2}  \Ex_\vartheta
\left|Z_T^{1/2}\left(u_2\right)-Z_T^{1/2}\left(u_1\right)\right|^2\leq
C\,\left(1+R^2\right).
\end{equation}
\end{lemma}
{\bf Proof.} Let us put $\vartheta _1=\vartheta +T^{-3/2}u_1,\vartheta
_2=\vartheta +T^{-3/2}u_2$ and $\delta \left(t,x\right)$ is defined below in
\eqref{a3}. Then the estimate \eqref{A1} with $m=1$ allows us to write
\begin{align*}
&\Ex_{\vartheta} \left|Z_T^{1/2}\left(u_2\right)-Z_T^{1/2}\left(u_1\right)\right|^{2}\\ 
&\quad \leq
C_1\Ex_{\vartheta _1} \left(\int_{0}^{T}V_t\,    \delta
\left(t,X_t\right)^2{\rm d}t\right)^{2} +C_2\Ex_{\vartheta _1} \int_{0}^{T}V_t^{2}\,\delta
\left(t,X_t\right)^2{\rm d}t\\ 
&\quad \leq
C_1T \int_{0}^{T}\Ex_{\vartheta _1}V_t^{2}\,    \delta
\left(t,X_t\right)^{4}{\rm d}t +C_2\int_{0}^{T}\Ex_{\vartheta _1}V_t^{2}\,\delta
\left(t,X_t\right)^{2}{\rm d}t\\ 
&\quad =
C_1T \int_{0}^{T}\Ex_{\vartheta _2}  \delta
\left(t,X_t\right)^{4}{\rm d}t +C_2\int_{0}^{T}\Ex_{\vartheta _2}\delta
\left(t,X_t\right)^{2}{\rm d}t.
\end{align*} 

As the derivative of the function $S\left(t\right)$ is bounded and we have
\eqref{15} we can write 
$$
T \int_{0}^{T}\Ex_{\vartheta _2}  \delta
\left(t,X_t\right)^{4}{\rm d}t\leq C T\left(u_2-u_1\right)^4 T^{-6}
\int_{0}^{T}t^4{\rm d}t\leq C\left|u_2-u_1\right|^4 
$$
and similary
$$
 \int_{0}^{T}\Ex_{\vartheta _2}  \delta
\left(t,X_t\right)^{2}{\rm d}t\leq C \left(u_2-u_1\right)^2 T^{-3}
\int_{0}^{T}t^2{\rm d}t\leq C\left|u_2-u_1\right|^2.
$$
Therefore this lemma is proved.

\begin{lemma}
\label{L4} For sufficiently large $T$ we have
\begin{equation}
\label{20}
\sup_{\vartheta \in \KK}\Ex_\vartheta Z_T^{1/2}\left(u\right)\leq e^{-\kappa\left| u\right|^{2/3}}
\end{equation}

\end{lemma}
{\bf Proof.} Let us put $\vartheta _2=\vartheta +T^{-3/2}u, \vartheta
_1=\vartheta $. We can write
\begin{align*}
&\Ex_\vartheta Z_T^{1/2}\left(u\right)\\
&\quad =\Ex_\vartheta
\exp\left\{\int_{0}^{T}{\frac{\delta \left(t,X_t\right)}{2}}\,{\rm d}W_t-\int_{0}^{T}{\frac{\delta
  \left(t,X_t\right)^2}{8}}\,{\rm d}t-\int_{0}^{T}{\frac{\delta
  \left(t,X_t\right)^2}{8}}\,{\rm d}t\right\} \\
&\quad \leq 
\exp\left\{-\frac{1}{8K}\int_{0}^{T}\left[S\left(\vartheta
  t+T^{-3/2}ut\right)-S\left(\vartheta t\right)\right]^2\,{\rm d}t\right\}  
\end{align*}
because
$$
\Ex_\vartheta
\exp\left\{\int_{0}^{T}{\frac{\delta \left(t,X_t\right)}{2}}\,{\rm d}W_t-\frac{1}{2}\int_{0}^{T}{\frac{\delta
  \left(t,X_t\right)^2}{4}}\,{\rm d}t\right\} =1
$$
and by condition  \eqref{15} we have as well 
$$
\int_{0}^{T}{{\delta
  \left(t,X_t\right)^2}}\,{\rm d}t\geq \frac{1}{K}\int_{0}^{T}\left[S\left(\vartheta
  t+T^{-3/2}ut\right)-S\left(\vartheta t\right)\right]^2\,{\rm d}t.
$$
For the last integral according to \eqref{inc} we have ($z=\vartheta ^{-1}T^{-1/2}u$)
\begin{align*}
\int_{0}^{T}\left[S\left(\vartheta
  t+T^{-3/2}ut\right)-S\left(\vartheta t\right)\right]^2\,{\rm d}t\geq
c\,T\,\frac{\frac{u^2}{\vartheta ^2T}}{1+\frac{u^2}{\vartheta ^2T }}.
\end{align*}
Further, if $u^2\leq \vartheta ^2T$, then 
$$
T\,\frac{\frac{u^2}{\vartheta ^2T}}{1+\frac{u^2}{\vartheta ^2T }}\geq
\frac{u^2}{2\vartheta ^2}, 
$$
and if $u^2> \vartheta ^2T$, then
$$
T\,\frac{\frac{u^2}{\vartheta ^2T}}{1+\frac{u^2}{\vartheta ^2T }}\geq
\frac{T}{2}\geq \frac{\left|u\right|^{2/3}}{2\left(\beta -\alpha \right)^{2/3}} 
$$
because $\left|u\right|\leq T^{3/2}\left(\beta -\alpha \right)$. 
Therefore
$$
\frac{1}{8K}\int_{0}^{T}{{\delta
  \left(t,X_t\right)^2}}\,{\rm d}t\geq \kappa \left|u\right|^{2/3}
$$
with some positive $\kappa $. 

The properties of the likelihood ratio process established in Lemmas 1-3 allow
us to apply the Theorems 3.1.1, 3.1.3 and 3.2.1 in \cite{IH81} and to obtain
all properties of the MLE and BE announced in the Theorem \ref{T1}. 

\subsection{Discontinuous trend}

We have the same model of observed periodic diffusion process \eqref{14} but
the function $S\left(t\right),t\geq 0$ is now  discontnuous. More precizely,
the following condition holds.

\bigskip

{\bf C.}  {\it The function $S\left(\cdot \right)$ is periodic with period 1,
  is continuously differentiable everywhere except the points $\tau_* +k$
  ($\tau _*\in \left(0,1\right), k=0,1,2,\ldots $) and at the points $\tau_*+k
  $ it has the  left
  and right limits $S\left(\tau_* -\right)$ and $S\left(\tau_* +\right)$
  respectively,} $S\left(\tau_* +\right)-S\left(\tau_* -\right)=r\not=0$.

\bigskip

The likelihood ratio random function $L\left(\vartheta ,X^T\right), \vartheta
\in \Theta $ is continuous with probability 1 (see Lemma \ref{L5} below),
hence the solution of equation  \eqref{MLE} exists and the both estimators
$\hat\vartheta _T$ and  $\tilde\vartheta _T$ are well defined. 

The limit distributions of these estimators are described with the help of the
random variables $\hat u$ and $\tilde u$ defined in \eqref{10} where
$Z\left(u\right)$  is given by \eqref{9} with
$$
\gamma^2 ={\left[S\left(\tau _*+\right)-S\left(\tau
  _*-\right)\right]^2}{}\int_{-\infty }^{\infty }\frac{f_\vartheta \left(\tau
 _*,x\right)}{2\,\sigma \left(x\right)^2}\,{\rm d}x .
$$

 The lower bound on the meansquare risk of all estimators is similar to \eqref{bound1}
:
\begin{equation*}
\label{bound2}
\Liminf_{\delta \rightarrow 0}\Liminf_{T \rightarrow \infty
}\sup_{\left|\vartheta -\vartheta _0\right|\leq \delta }T^4 \Ex_\vartheta
\left|\bar\vartheta _T-\vartheta \right|^2\geq \Ex_{\vartheta _0}\tilde u^2.
\end{equation*}
For the proof of more general result  see \cite{IH81}, Section 1.9. In our case
this inequality can be prooved in three lines if we suppose that we have  already
proved the uniform convergence of moments of the BE 
$$
T^4 \Ex_\vartheta
\left|\tilde\vartheta _T-\vartheta \right|^2\longrightarrow \Ex_{\vartheta}
\left|\tilde u\right|^2
$$
(Theorem \ref{T2} below)
as follows. Let us denote $p_\delta \left(\theta \right), \theta \in
\left[\vartheta _0-\delta ,\vartheta _0+\delta \right]$ a density function and
$\tilde \theta _T$ the corresponding bayesian estimator. Then we can write
\begin{align*}
&\sup_{\left|\vartheta -\vartheta _0\right|\leq \delta }T^4 \Ex_\vartheta
\left|\bar\vartheta _T-\vartheta \right|^2\geq T^4 \int_{\vartheta
  _0-\delta}^{\vartheta _0+\delta}\Ex_\theta 
\left|\bar\vartheta _T-\theta \right|^2 p_\delta \left(\theta \right){\rm
  d}\theta \\
&\quad \geq T^4 \int_{\vartheta
  _0-\delta}^{\vartheta _0+\delta}\Ex_\theta 
\left|\tilde\theta _T-\theta \right|^2 p_\delta \left(\theta \right){\rm
  d}\theta\longrightarrow \int_{\vartheta
  _0-\delta}^{\vartheta _0+\delta}\Ex_\theta 
\left|\tilde u\right|^2 p_\delta \left(\theta \right){\rm
  d}\theta=\Ex_{\vartheta_0 }
\left|\tilde u\right|^2.
\end{align*}

This bound allows us to  call an estimator $\vartheta _T^\star$
 asymptotically efficient if for all $\vartheta 
_0\in \Theta $ we have the equality
\begin{equation}
\label{asefd}
\lim_{\delta \rightarrow 0}\lim_{T \rightarrow \infty
}\sup_{\left|\vartheta -\vartheta _0\right|\leq \delta }T^4 \Ex_\vartheta
\left|\vartheta _T^\star-\vartheta \right|^2= \Ex_{\vartheta _0}\left|\tilde u\right|^2.
\end{equation}

Remind that the last value is known : $\Ex_{\vartheta }\left|\tilde
u\right|^2\approx 19,3\,\gamma ^{-4} $ \cite{RS} and is less than the
similar quantity for the MLE $\Ex_{\vartheta}\left|\hat
u\right|^2=26\,\gamma ^{-4}$ \cite{Ter}.

\begin{theorem}
\label{T2} 
Let the conditions {\bf A} and  {\bf C} be fulfilled. Then  the MLE $\hat \vartheta _T$ and
BE $\tilde\vartheta _T$ have the following properties uniformly on compacts
$\KK\subset \Theta $ :
\begin{itemize}
\item 
The both estimators are consistent.
\item  They have different limit distributions
\begin{equation*}
T^{2}\left(\hat\vartheta _T-\vartheta \right)\Longrightarrow \hat u,\qquad
T^{2}\left(\tilde\vartheta _T-\vartheta \right)\Longrightarrow \tilde u.
\end{equation*}
\item The convergence of moments take place : for any $p>0$
$$
\lim_{T\rightarrow \infty }T^{2p}\;\Ex_\vartheta \left|\hat\vartheta
_T-\vartheta \right|^p=\Ex_\vartheta \left|\hat u\right|^p,\qquad \lim_{T\rightarrow
  \infty }T^{2p}\;\Ex_\vartheta \left|\tilde\vartheta 
_T-\vartheta \right|^p=\Ex_\vartheta \left|\tilde u\right|^p,
$$
\item BE are asymptotically efficient in the sense \eqref{asefd}. 
\end{itemize}
\end{theorem}

{\bf Proof.}  The normalized likelihood ratio is now
\begin{align*}
Z_T\left(u\right)=\frac{L\left(\vartheta +T^{-2}u,X^T\right)}{L\left(\vartheta
  ,X^T\right)},\qquad u\in 
\UU_T=\left(T^{2}\left(\alpha -\vartheta 
\right),T^{2}\left(\beta  -\vartheta \right) \right).
\end{align*} 
We have the following result (see \cite{HK10c}, Theorem 1.2)
\begin{lemma}
\label{L5a} The finite dimensional distributions of the random process
$Z_T\left(u\right)$ converge to the finite dimensional distributions of the
process $Z\left(u\right)$ uniformly on compacts $\KK\in \Theta $.
\end{lemma}
To explain this convergence we can write the log-likelihood ratio as follows:
  (below $\vartheta _u=\vartheta +T^{-2}u, u>0$)
\begin{align*}
&\ln Z_T\left(u\right)=\int_{0}^{T}\frac{S\left(\vartheta
  _ut\right)-S\left(\vartheta t\right)}{\sigma \left(X_t\right)}{\rm
  d}W_t-\frac{1}{2}\int_{0}^{T}\left(\frac{S\left(\vartheta
  _ut\right)-S\left(\vartheta t\right)}{\sigma \left(X_t\right)}\right)^2{\rm
  d}t\\
&\, =\sum_{k=0}^{\left[T\vartheta \right]}\int_{\frac{\tau
      _*+k}{\vartheta_u}}^{\frac{\tau _*+k}{\vartheta }
  }\frac{S\left(\vartheta 
  _ut\right)-S\left(\vartheta t\right)}{\sigma \left(X_t\right)}{\rm
  d}W_t-\sum_{k=0}^{\left[T\vartheta \right]}\int_{\frac{\tau _*+k}{\vartheta
    _u}}^{\frac{\tau _*+k}{\vartheta } }\frac{\left(S\left(\vartheta
  _ut\right)-S\left(\vartheta t\right)\right)^2}{2\,\sigma \left(X_t\right)^2}{\rm
  d}t+o\left(1\right)\\
&\, =\sum_{k=0}^{\left[T\vartheta \right]}\left[r\int_{\frac{\tau
      _*+k}{\vartheta_u}}^{\frac{\tau _*+k}{\vartheta }  }\frac{{\rm
      d}W_t}{\sigma \left(X_t\right)} -\frac{r^2}{2}\int_{\frac{\tau _*+k}{\vartheta
    _u}}^{\frac{\tau _*+k}{\vartheta } }\frac{{\rm d}t}{\sigma
    \left(X_t\right)^2}\right]+o\left(1\right)\Longrightarrow  \gamma
  \,W\left(u\right)-\frac{\gamma ^2}{2} \,\left|u\right|
\end{align*}
because 
\begin{align*}
&\sum_{k=0}^{\left[T\vartheta \right]}\int_{\frac{\tau _*+k}{\vartheta
    _u}}^{\frac{\tau _*+k}{\vartheta } }\frac{{\rm d}t}{\sigma 
    \left(X_t\right)^2}=\frac{u}{\vartheta ^2T^2}\sum_{k=0}^{\left[T\vartheta
    \right]}\;\frac{k}{\sigma \left(X_{\frac{\tau _*+k}{\vartheta
  }}\right)^2}+o\left(1\right) \\
&\quad \longrightarrow \frac{u}{2}\int_{-\infty }^{\infty }\frac{f_\vartheta
    \left(\frac{\tau _*}{\vartheta },x\right)}{\sigma \left(x\right)^2}\;{\rm
    d}x =\frac{u}{2}\int_{-\infty }^{\infty }\frac{f_1
    \left({\tau _*},x\right)}{\sigma \left(x\right)^2}\;{\rm
    d}x
\end{align*}
Here $\left[T\vartheta \right]$ is the integer part of $T\vartheta $. The
asymptotic normality of the stochastic integral follows from this 
convergence (central limit theorem). For the details see \cite{HK10c}.

\begin{lemma}
\label{L5} The following inequaulity holds
\begin{equation}
\label{28}
\sup_{\vartheta \in\KK}\Ex_\vartheta
\left|Z_T^{1/4}\left(u_2\right)-Z_T^{1/4}\left(u_1\right)\right|^4\leq
C\,\left|u_2-u_1\right|^2 
\end{equation}
\end{lemma}
{\bf Proof.} According to \eqref{A1} with $m=2$ we have
\begin{align*}
&\Ex_\vartheta
\left|Z_T^{1/4}\left(u_2\right)-Z_T^{1/4}\left(u_1\right)\right|^4\\
&\qquad  \leq
C_1\Ex_{\vartheta _1} \left(\int_{0}^{T}V_t\,\delta 
\left(t,X_t\right)^2{\rm d}t\right)^{4} +C_2\Ex_{\vartheta _1} \left(\int_{0}^{T}V_t^2\,\delta
\left(t,X_t\right)^2{\rm d}t\right)^{2}.
\end{align*}

Let us put $\vartheta _{u_1}=\vartheta +T^{-2}u_1,\vartheta _{u_2}=\vartheta
+T^{-2}u_2 $, $N=\left[T\vartheta
    _{u_1}\right]$ and consider the case $ 0<u_1<u_2$. Then we can write
\begin{align*}
\int_{0}^{T}V_t\,\delta \left(t,X_t\right)^2{\rm
  d}t=\sum_{k=0}^{N-1}\left[\int_{\frac{k}{\vartheta _{u_1}}}^{\frac{\tau
      _*+k}{\vartheta _{u_2}}}+\int_{\frac{\tau _*+k}{\vartheta
      _{u_2}}}^{\frac{\tau _*+k}{\vartheta _{u_1}}}+\int_{\frac{\tau
      _*+k}{\vartheta _{u_1}}}^{\frac{k+1}{\vartheta
      _{u_1}}}\right]V_t\,\delta 
\left(t,X_t\right)^2{\rm d}t.
\end{align*}
The function $\delta \left(t,X_t\right)^2$ on the intervals 
$$
\left[  \frac{k}{\vartheta _{u_1}},\frac{\tau
      _*+k}{\vartheta _{u_2}}\right]\qquad {\rm and}\qquad \left[\frac{\tau
      _*+k}{\vartheta _{u_1}},\frac{k+1}{\vartheta
      _{u_1}}\right]
$$
is continuously differentiable on $\vartheta $ and therefore is majorated as
follows
$$
\delta \left(t,X_t\right)^2\leq C\,t^2\;\frac{\left(u_2-u_1\right)^2}{T^4}.
$$
Further, we have
$$
\int_{\frac{\tau _*+k}{\vartheta      _{u_2}}}^{\frac{\tau _*+k}{\vartheta _{u_1}}}V_t\,\delta 
\left(t,X_t\right)^2{\rm d}t\leq C\;\left[S\left(\tau _*+\right)-S\left(\tau
  _*-\right)\right]^2  \int_{\frac{\tau _*+k}{\vartheta
    _{u_2}}}^{\frac{\tau _*+k}{\vartheta _{u_1}}}V_t\,{\rm d}t .
$$
Hence
\begin{align*}
&\Ex_{\vartheta _1}\left(\sum_{k=0}^{N-1}\int_{\frac{\tau _*+k}{\vartheta
    _{u_2}}}^{\frac{\tau _*+k}{\vartheta _{u_1}}}V_t\,{\rm d}t\right)^4\leq
C\,N^3\sum_{k=0}^{N-1}\Ex_{\vartheta _1}\left(\int_{\frac{\tau _*+k}{\vartheta 
    _{u_2}}}^{\frac{\tau _*+k}{\vartheta _{u_1}}}V_t\,{\rm d}t\right)^4 \\
&\qquad \leq  C\,N^3\sum_{k=0}^{N-1} \frac{\left(\tau
  _*+k\right)^3\left(u_2-u_1\right)^3}{T^6}   \int_{\frac{\tau _*+k}{\vartheta  
    _{u_2}}}^{\frac{\tau _*+k}{\vartheta _{u_1}}}\Ex_{\vartheta
  _1}V_t^4\,{\rm d}t\\
&\qquad \leq  C\,N^3\sum_{k=0}^{N-1} \frac{\left(\tau
  _*+k\right)^4\left(u_2-u_1\right)^4}{T^8}\leq C\,\left(u_2-u_1\right)^4.
\end{align*}
Remind that $\Ex_{\vartheta   _1}V_t^4=1$. For the second integral the similar
arguments provide 
$$
\Ex_{\vartheta _1} \left(\int_{0}^{T}V_t^2\,\delta
\left(t,X_t\right)^2{\rm d}t\right)^{2}\leq C\,\left(u_2-u_1\right)^2.
$$
Now \eqref{28} follows from the last two estimates.

\begin{lemma}
\label{L7} For sufficiently large $T$ we have
\begin{equation}
\label{25}
\sup_{\vartheta \in \KK}\Ex_\vartheta Z_T^{1/2}\left(u\right)\leq e^{-\kappa\left| u\right|^{1/2}}
\end{equation}

\end{lemma}
{\bf Proof.} Following  the proof of  Lemma \ref{L4} (with $\vartheta
_u=\vartheta +T^{-2}u$) we obtain  the estimates 
\begin{align*}
&\Ex_\vartheta Z_T^{1/2}\left(u\right) \leq 
\exp\left\{-\frac{1}{8K}\int_{0}^{T}\left[S\left(\vartheta
  t+T^{-2}ut\right)-S\left(\vartheta t\right)\right]^2\,{\rm d}t\right\}  
\end{align*}
and 
$$
\int_{0}^{T}{{\delta
  \left(t,X_t\right)^2}}\,{\rm d}t\geq \frac{1}{K}\int_{0}^{T}\left[S\left(\vartheta
  t+T^{-2}ut\right)-S\left(\vartheta t\right)\right]^2\,{\rm d}t.
$$
Further,  the estimate  \eqref{ind} allows us to write ($z=\vartheta ^{-1}T^{-1}u$)
\begin{align*}
\int_{0}^{T}\left[S\left(\vartheta
  t+T^{-2}ut\right)-S\left(\vartheta t\right)\right]^2\,{\rm d}t\geq
c\,T\,\frac{\frac{\left|u\right|}{\vartheta\,T}}{1+\frac{\left|u\right|}{\vartheta \,T }}.
\end{align*}
 If $\left|u\right|\leq \vartheta T$, then 
$$
T\,\frac{\frac{\left|u\right|}{\vartheta\,T}}{1+\frac{\left|u\right|}{\vartheta \,T }}\geq
\frac{\left|u\right|}{2\vartheta }, 
$$
and if $\left|u\right|> \vartheta T$, then
$$
T\,\frac{\frac{\left|u\right|}{\vartheta\,T}}{1+\frac{\left|u\right|}{\vartheta \,T }}\geq
\frac{T}{2}\geq \frac{\left|u\right|^{1/2}}{2\left(\beta -\alpha \right)^{1/2}} 
$$
because $\left|u\right|\leq T^{2}\left(\beta -\alpha \right)$. 
Therefore
$$
\frac{1}{8K}\int_{0}^{T}{{\delta
  \left(t,X_t\right)^2}}\,{\rm d}t\geq \kappa \left|u\right|^{1/2}
$$
with some positive $\kappa $. 

\bigskip

The convergence of the finite-dimensional distributions  of
the random function $Z_T\left(\cdot \right)$  together with \eqref{28} and \eqref{25}
allow us to cite the Theorems 1.10.1, 1.10.2, where the mentioned in the
Theorem \ref{T2} properties of estimators are proved.

\section{Auxiliary results}

{\bf Two lemmae.}
We remind here one estimate for the increaments of the likelihood ratio and
two lemmae which allowed us  to prove the Lemma 2,3,5,6. 

Let us introduce three diffusion processes
$$
{\rm d}X_t=\left[S\left(\vartheta _it\right)+b\left(X_t\right)\right]\,{\rm
  d}t+\sigma \left(X_t\right)\;{\rm 
  d}W_t,\quad X_0, 0\leq t\leq T,\quad i=0,1,2,
$$
and denote by 
 ${\bf P}_{\vartheta _i}^{\left(T\right)},i=0,1,2, $  the corresponding
measures induced by these processes in $\left({\cal C}\left[0,T\right],{\scr B}\left[0,T\right]\right)$. 
The Radon-Nikodym derivatives are denoted as 
$$
Z_i=\frac{{\rm d}{\bf P}_{\vartheta _i}^{\left(T\right)}}{{\rm d}{\bf
    P}_{\vartheta _0}^{\left(T\right)}}\left(X^T\right),\quad i=1,2,\quad
V_t=\left(\frac{{\rm d}{\bf P}_{\vartheta _2}^{\left(t\right)} }{{\rm d}{\bf
    P}_{\vartheta _1}^{\left(t\right)}}\left(X^t\right)\right)^{1/2m}  
$$
where $m\geq 1$ is some integer. 
Below we  put 
\begin{equation}
\label{a3}
\delta \left(t,x\right)= \frac{S\left(\vartheta _2t\right)- S\left(\vartheta _1t\right)}{ \sigma
\left(x\right)}. 
\end{equation}
Remind  that we suppose \eqref{15} and that the function
$S\left(t\right),t\geq 0$ is bounded, hence the function $\delta
\left(t,x\right)$ is bounded too and by Lemma 1.13 in \cite{Kut04} we have
the following result.

{\it  There exist constants $C_1\left(m\right)$ and $C_2\left(m\right)$ such
that}
\begin{align}
\label{A1}
\Ex_{\vartheta _0} \left|Z_2^{1/2m}-Z_1^{1/2m}\right|^{2m}&\leq
C_1\left(m\right)\Ex_{\vartheta _1} \left(\int_{0}^{T}V_t\,\delta
\left(t,X_t\right)^2{\rm d}t\right)^{2m}\nonumber\\ 
&\qquad +C_2\left(m\right)\Ex_{\vartheta _1} \left(\int_{0}^{T}V_t^2\,\delta
\left(t,X_t\right)^2{\rm d}t\right)^{m}
\end{align}

\bigskip

The exponential decreasing of the tails of $Z_T\left(u\right)$ are verified
with the help of the following two lemmas.

\begin{lemma} {\rm (Ibragimov and Khasminskii)} \label{LA1}
Let $h\left(t\right)$ be a nonconstant continuously differentiable periodic
function. Then for all $T$ sufficiently large and for some constant
$c >0$, the inequality
\begin{equation}
\label{inc}
\frac{1}{T}\int_{0}^{T}\left[h\left(t+\frac{z}{T}t\right)-h\left(t\right)\right]^2{\rm
  d}t\geq c\;\frac{z^2}{1+z^2}
\end{equation}
is valid.
\end{lemma}
For the proof see \cite{IH81}, Lemma 3.5.3. 

\bigskip

The similar result for discontinuous function is given in the following
Lemma. 

\begin{lemma}  \label{LA2} 
Let $S\left(t\right)$ satisfies the condition {\bf B}. Then for all $T$
sufficiently large and for some constant 
$c >0$, the inequality
\begin{equation}
\label{ind}
\frac{1}{T}\int_{0}^{T}\left[S\left(t+\frac{z}{T}t\right)-S\left(t\right)\right]^2{\rm
  d}t\geq c\;\frac{\left|z\right|}{1+\left|z\right|}
\end{equation}
is valid.
\end{lemma}
The proof of this lemma is a modification of the proof of lemma \ref{LA1},
which can be found in \cite{Kut98}, Lemma 5.7. 
 
\bigskip

{\bf On uniform convergence}. The uniform in $\vartheta \in \Theta $
convergence (condition {\bf A4}) 
$$
\frac{1}{T}\int_{0}^{T}h\left(\vartheta ,t,X_t\right)\;{\rm d}t\longrightarrow \frac{1}{\tau }
\int_{-\infty }^{\infty }\int_{0}^{\tau }h\left(\vartheta ,t,x\right)
\;f_\vartheta \left(t,x\right){\rm d}t\;{\rm d}x\equiv A\left(\vartheta \right).
$$
means, that for any $\varepsilon >0$ we have
\begin{align}
\label{cc}
\sup_{\vartheta \in \Theta }\Pb_\vartheta \left\{
\left|\frac{1}{T}\int_{0}^{T}\left[h\left(\vartheta
  ,t,X_t\right)-A\left(\vartheta \right)\right]\;{\rm 
  d}t \right|>\varepsilon  \right\}\longrightarrow 0
. \end{align}
Let us denote by  $H\left(\vartheta ,t,x\right)$ the solution of the following equation
\begin{align*}
\frac{\partial H}{\partial
  t}+\left[S\left(\vartheta t\right)+b\left(x\right)\right]\frac{ \partial H}{\partial
  x} +\frac{\sigma \left(x\right)^2}{2}\frac{ \partial^2 H}{\partial
  x^2}=h\left(\vartheta ,t,x\right)-A\left(\vartheta \right)
\end{align*}
Then we can write
\begin{align*}
\frac{1}{T}&\int_{0}^{T}\left[h\left(\vartheta ,t,X_t\right)-A\left(\vartheta \right)\right]\;{\rm
  d}t=\frac{H\left(\vartheta ,T,X_T\right)-H\left(\vartheta ,0,X_0\right)}{T}\\
&-\frac{1}{T}\int_{0}^{T}H'_x\left(\vartheta ,t,X_t\right)\sigma
\left(X_t\right) {\rm d}W_t
\end{align*}

Hence
\begin{align*}
&\Ex_\vartheta \left(\frac{1}{T}\int_{0}^{T}\left[h\left(\vartheta
    ,t,X_t\right)-A\left(\vartheta \right)\right]\;{\rm d}t\right)^2 \leq
  2\frac{\Ex_\vartheta \left(H\left(\vartheta ,T,X_T\right)-H\left(\vartheta ,0,X_0\right)\right)^2
  }{T^2}\\ &\qquad \qquad
  +\frac{2}{T^2}\int_{0}^{T}\Ex_\vartheta\left(H'_x\left(\vartheta ,t,X_t\right)\sigma
  \left(X_t\right) \right)^2{\rm d}t
\end{align*}
It is sufficient to suppose that the last expectations are bounded uniformly
in $\vartheta \in\Theta $ and apply in \eqref{cc} the  Tchebyshev inequality.

\section{Discussion}

{\bf Choice of the signal.} Let us consider the equation \eqref{7} with two types of modulation : p) phase
$S\left(\vartheta ,t\right)=S\left(t-\vartheta \right)$ and f) frequency
$S\left(\vartheta ,t\right)=S\left(\vartheta t\right)$ in two situations :
smooth and discontinuous. Then from the results obtained in
\cite{HK10b}-\cite{HK10c} and presented in this work it follows that we have
four problems with four different rates
\begin{align*}
&{\rm smooth}\qquad \;\;\quad \quad (p) \quad \Ex_\vartheta \left(\hat\vartheta _T-\vartheta
\right)^2 \sim \frac{c}{T}, \quad  \quad (f) \quad \Ex_\vartheta \left(\hat\vartheta _T-\vartheta
\right)^2 \sim \frac{c}{T^3}\\ 
&{\rm discontinuous}\quad (p)\quad \Ex_\vartheta \left(\hat\vartheta _T-\vartheta
\right)^2 \sim \frac{c}{T^2}, \quad  \qquad (f)\quad \Ex_\vartheta \left(\hat\vartheta _T-\vartheta
\right)^2 \sim \frac{c}{T^4}.
\end{align*}
It is natural to ask: {\it  how far can we go in the rate of convergence?  What is
the best choice of the signal and what is the best rate?} 

The similar statement for the {\it signal in white Gaussian noise} problem
$$
{\rm d}X_t=S\left(\vartheta ,t\right){\rm d}t+{\rm d}W_t,\quad X_0=0,\quad  0\leq
t\leq T,\quad  \vartheta \in \left[0,1\right]
$$
was considered by M. Burnashev \cite{Bur4}. It was shown that for signals
satisfying
\begin{equation}
\label{r}
\frac{1}{T}\int_{0}^{T}S\left(\vartheta ,t\right)^2{\rm d}t\leq L
\end{equation}
the best choice yields ($T\rightarrow \infty $)
$$
\inf_{S,\bar\vartheta _T}\sup_{\vartheta \in \left[0,1\right]}\Ex_\vartheta \left(\hat\vartheta _T-\vartheta
\right)^2 =\exp\left\{-\frac{L}{6}T\left(1+o\left(1\right)\right)\right\}.
$$
Therefore the rate can be even exponential. The similar result was obtained
for inhomogenneous Poisson processes too \cite{BK}. It follows that if the
diffusion coefficient $\sigma \left(x\right)^2\equiv 1$ and the signal
$S\left(\vartheta ,t\right)$ in the equation \eqref{7} satisfies the condition
\eqref{r}, then we have the same result with exponential rate.

\bigskip

{\bf Generalisations.} There are several generalisations which can be done by
the direct calculations similar to one given above.

If the function $S\left(t\right)$ has jumps in $k$ points $0<\tau
_1,\ldots,\tau _k<\tau  $ and is continuously differentiable between these
points, then the estimators have the same   asymptotic properties as described
in the Theorem \ref{T2} but the constant is
$$
\gamma ^2=\sum_{l=1}^{k}\int_{}^{}\frac{\left[S\left(\tau _l+\right)-S\left(\tau
  _{l}-\right)\right]^2}{\sigma \left(x \right)^2}\,f_\vartheta \left(\tau
_l+\vartheta ,x\right)\,{\rm d}x .
$$
The problem becames a bit more complicate if
$$
{\rm d}X_t=\sum_{l=1}^{k}S_l\left(t-\vartheta _l\right){\rm
  d}t+b\left(X_t\right){\rm d}t+ \sigma \left(X_t\right){\rm d}W_t
$$ and $\vartheta =\left(\vartheta _1,\ldots,\vartheta _k\right)$ but can be
done too. The limit likelihood ratio is a product of $k$ one-dimensional
likelihood ratios.  See the details in \cite{Kut98}, where the similar
problems were considered for periodic Poisson processes.

 \bigskip

{\bf Three-dimensional parameter.} In both problems studied above we supposed that the unknown parameter is
one-dimensional. It is interesting to see the properties of estimators, say,
in smooth case when the we have three dimensional parameter $\vartheta
=\left(\rho ,\omega ,\varphi \right)$ in the model of observations
$$
{\rm d}X_t=\rho \,\sin\left(2\pi \omega \,t+\varphi \right)\,{\rm d}t+b(X_t){\rm
  d}t+\sigma \left(X_t\right){\rm d}W_t,\quad X_0=0, 0\leq t\leq T,
$$
 i.e.: we have to estimate the 
amplitude $\rho  $, frequency $\omega $ and phase $\varphi $ of the signal
$S\left(\vartheta ,t\right)=\rho \,\sin\left(\omega \,t+\varphi \right) $.

The functions $b\left(\cdot \right)$ and $\sigma \left(\cdot \right)$ satisfy
the condition {\bf A}. This problem of parameter estimation is regular and the
technique developped in this work together with calculus presented in Example 4
Section 3.5 in  \cite{IH81} allows us to show that the MLE $\hat\vartheta _T$ is
consistent and asymptotically normal :
$$
\sqrt{T}\left(\hat\rho _T-\rho \right)\Longrightarrow \eta ,\quad
     {T^{3/2}}\left(\hat\omega  _T-\omega  \right)\Longrightarrow \xi ,\quad
     \sqrt{T}\left(\hat\varphi  _T-\varphi  \right)\Longrightarrow \zeta  
$$

\vskip1.5cm
\small
Reinhard H\"opfner\\
Institut f\"ur Mathematik, Universit\"at Mainz, D--55099 Mainz\\ 
{\tt hoepfner@mathematik.uni-mainz.de}\\
{\tt http://www.mathematik.uni-mainz.de/$\sim$hoepfner}

\vskip0.5cm
 Yury A.\ Kutoyants\\
Laboratoire de Statistique et Processus, Universit\'e du Maine, F--72085 Le Mans\\ 
{\tt kutoyants@univ-lemans.fr}\\
{\tt http://lmm.univ-lemans.fr/spip.php?article22}
 
\end{document}